\documentclass[12pt]{scrartcl}

\usepackage[english]{babel}
\RequirePackage[T1]{fontenc}
\RequirePackage[utf8]{inputenc}

\RequirePackage{geometry}
\RequirePackage{palatino}
\RequirePackage[osf,sc]{mathpazo} 
\RequirePackage{mathrsfs} 
\RequirePackage{stmaryrd} 
\RequirePackage{microtype} 

\RequirePackage[leqno]{amsmath}
\RequirePackage{amsfonts}
\RequirePackage{amssymb}
\usepackage[svgnames]{xcolor} 
\RequirePackage[pdftex]{graphicx}

\RequirePackage{marginnote}  
\usepackage[alwaysadjust]{paralist}
\setdefaultenum{\footnotesize (1.)}{\footnotesize (a)}{\footnotesize (i.)}{\footnotesize A.}

\RequirePackage[pdftex,pagebackref,plainpages=false, pdfpagelabels,unicode]{hyperref}

\newenvironment{prose}{\begin{quotation}\footnotesize\noindent[[}
{]]\end{quotation}}

\providecommand{\Marginnote}[2][0pt]{\marginnote{\small {\color{RoyalBlue}#2}}[#1]}

\RequirePackage[amsmath,thmmarks,thref,hyperref]{ntheorem}

{
\theoremstyle{break}
\newtheorem{preTheorem}{Theorem}[section]
\newtheorem{preLemma}[preTheorem]{Lemma}
\newtheorem{preCorollary}[preTheorem]{Corollary}
\newtheorem{preProposition}[preTheorem]{Proposition}
}
{\theoremstyle{plain}
\theorembodyfont{\normalfont}
\newtheorem{preDefinition}[preTheorem]{Definition}
\newtheorem{preRemark}[preTheorem]{Remark}
\newtheorem{preExample}[preTheorem]{Example}
\newtheorem{preQuestion}[preTheorem]{Question}
}

\newenvironment{Theorem}{\bgroup \begin{preTheorem} }{\end{preTheorem} \egroup}
\newenvironment{Lemma}{\bgroup \begin{preLemma} }{\end{preLemma} \egroup}
\newenvironment{Corollary}{\bgroup \begin{preCorollary} }{\end{preCorollary} \egroup}
\newenvironment{Proposition}{\bgroup \begin{preProposition}}{\end{preProposition} \egroup}

\newenvironment{Definition}{\bgroup \begin{preDefinition} }{\end{preDefinition} \egroup}

\newenvironment{Question}{\bgroup \begin{preQuestion} }{\end{preQuestion} \egroup}

\newcommand{\plist}{\begin{compactenum}}
\newcommand{\pliste}{\end{compactenum}}

{
\theoremstyle{nonumberplain}
\theoremheaderfont{\normalfont\slshape}
\theorembodyfont{\upshape}
\theoremseparator{.}
\theoremsymbol{\ensuremath{_\Box}}
\newtheorem{proof}{Proof}
}

{
\theoremstyle{nonumberplain}
\theoremheaderfont{\normalfont}
\theorembodyfont{\upshape \footnotesize}
\theoremseparator{}
\theoremsymbol{ \ensuremath{_\boxtimes}}

\newtheorem{spoiler}{Summary.}
}

\newcommand{\hl}[1]{{\color{RoyalBlue}#1}}

\mathchardef\mathhyphen="2D 

\newcommand{\gflash}{{\color{RoyalBlue}\ensuremath{\lightning}}}

\newcommand{\widebar}[1]{\overline{#1}}
\newcommand{\bbar}{\widebar} 

\newcommand{\matn}{\mathbb{N}}
\newcommand{\matz}{\mathbb{Z}}

\newcommand{\matf}{\mathbb{F}}

\newcommand{\mathp}{\mathfrak{P}}

\newcommand{\restr}{\!\upharpoonright\!}
\newcommand{\notnil}{\not= \emptyset}
\newcommand{\nil}{\emptyset}
\newcommand{\phii}{\varphi}

\newcommand{\mal}{\cdot}

\newcommand{\bsupp}{\mathrm{bsupp}}
\newcommand{\bmax}{\mathrm{bmax}}
\newcommand{\bmin}{\mathrm{bmin}}

\newcommand{\xsupp}{\mathrm{\mathhyphen supp}}
\newcommand{\xmin}{\mathrm{\mathhyphen min}}
\newcommand{\xmax}{\mathrm{\mathhyphen max}}
\newcommand{\ZFC}{\textsl{ZFC}}

\newcommand{\dotcup}{\dot{\cup}}

\hypersetup{
   pdftitle={On strongly summable ultrafilters},    
   pdfauthor={Peter Krautzberger},     
   pdfsubject={Some new results on strongly summable ultrafilters.},   
   pdfkeywords={ mathematics, ultrafilters, stone-cech compactification, combinatorics, set theory, semigroup, idempotent ultrafilter,union ultrafilter, summable ultrafilter, strongly summable ultrafilter}, 
   linkbordercolor = RoyalBlue,
   citebordercolor = RoyalBlue,
   urlbordercolor = RoyalBlue
    }

\title{On strongly summable ultrafilters}
\author{Peter Krautzberger\thanks{Mathematics Department, University of Michigan, Ann Arbor, pkrautzb@umich.edu, Partially supported by DFG-grant KR 3818; Subject classification 03E75 (Primary) 54D80, 05D10 (Secondary)} }
\date{\today}

\begin{document}

\maketitle
\begin{abstract}
We present some new results on strongly summable ultrafilters. As the main result, we extend a theorem by N.~Hindman and D.~Strauss on writing strongly summable ultrafilters as sums.
\end{abstract}

\emph{This preprint has been accepted for publication by the New York Journal of Mathematics and is available at \url{http://nyjm.albany.edu/j/2010/Vol16.htm}.}

\section*{Introduction}

The equivalent notions of strongly summable and union ultrafilters have been important examples of idempotent ultrafilters ever since they were first conceived in \cite{Hindman72}, \cite{Blass87-1} respectively. Their unique properties have been applied in set theory, algebra in the Stone-\v Cech compactification and set theoretic topology. For example, strongly summable ultrafilters were, in a manner of speaking, the first idempotent ultrafilters known, cf.~\cite{Hindman72} and \cite[notes to  Chapter 5]{HindmanStrauss}; they were the first strongly right maximal idempotents known and they are the only known class of idempotents with a maximal group isomorphic to $\matz$. Their existence is independent of \ZFC, since it implies the existence of (rapid) $P$-points, cf.~\cite{BlassHindman87}.\footnote{Also, as a strongly right maximal idempotent ultrafilter, the orbit closure of a strongly summable ultrafilter is an interesting example of a van Douwen space, cf.~\cite{HindmanStrauss02}; this will, however, not be relevant in what follows.}

The first part of this paper will focus on union ultrafilters for which we prove a new property; in the second part, this property is applied to strengthen a theorem on writing strongly summable ultrafilters as sums due to N.~Hindman and D.~Strauss \cite{HindmanStrauss95}, \cite[Chapter 12]{HindmanStrauss}.

The presentation of the proofs is inspired by \cite{Leron83} and \cite{Lamport95} splitting the proofs into different levels, at times adding {\footnotesize [[in the elevator]]} comments in between. The typesetting incorporates ideas from \cite{Tufte05} highlighting details in the proofs and structural remarks in the margin. Online discussion is possible through the author's website at \url{http://peter.krautzberger.info/papers}.

\section{Preliminaries}

Let us begin by giving a non-exhaustive selection of standard terminology in which we follow N.~Hindman and D.~Strauss \cite{HindmanStrauss}; for standard set theoretic notation we refer
to T.~Jech \cite{Jech}, e.g., natural numbers are considered as ordinals, i.e., $n = \{ 0, \ldots , n-1\}$. We work in \ZFC\ throughout. The main objects of this paper are \emph{(ultra)filters} on an infinite set $S$, i.e., (maximal) proper subsets of the power set $\mathp(S)$ closed under taking finite intersections and supersets.  
$S$ carries the discrete topology in which case the set of ultrafilters is $\beta S$, its Stone-\v Cech compactification. The Stone topology on $\beta S$ is generated by basic clopen sets induced by subsets $A\subseteq S$ in the form $ \bbar{A} := \{ p\in \beta S\ |\ A\in p\}$.
Filters are usually denoted by upper case Roman letters, mostly $F,G,H$, ultrafilters by lower case Roman letters, mostly $p,q,r,u$.

The set $S$ is always assumed to be the domain of a \Marginnote{(Partial) Semigroup}
\emph{(partial) semigroup} $(S,\cdot)$, i.e., the (partial) operation $\mal$ fulfills the associativity law  $ s\cdot (t \cdot v) = (s\cdot t) \cdot v $ (in the sense that if one side is defined, then so is the other and they are equal). For a partial semigroup $S$ and $s\in S$ the set of elements compatible with $s$ is denoted by $\sigma(s) := \{ t\in S\ |\ s\cdot t \mbox{ is defined} \}$.
A partial semigroup is also assumed to be \emph{adequate}, i.e., $\{ \sigma(s) \ |\ s\in S\}$ has the finite intersection property.
We denote the generated filter by $\sigma(S)$ and the corresponding closed subset of $\beta S$ by \index{$\delta S$}\index{semigroup!$\delta S$} $\delta S$.
For partial semigroups $S,T$ a map $\phii: S \rightarrow T$ is a \emph{partial semigroup homomorphism}\index{homomorphism!partial semigroup} if $\phii[\sigma(s)] \subseteq \sigma(\phii(s))$ and
\[
(\forall s\in S)(\forall s' \in \sigma(s))\ \phii(s\mal s') = \phii(s) \mal \phii(s').
\]
To simplify notation in a partial semigroup, $s \cdot t$ is always meant to imply $t \in \sigma(s)$. For $s\in S$, the restricted multiplication to $s$ from the left (right) is denoted by $\lambda_s$ ($\rho_s$).

It is easy to see that the operation of a partial semigroup can always be extended to a full semigroup operation by adjoining a (multiplicative) zero which takes the value of all undefined products. One key advantage of partial semigroups is that partial subsemigroups are usually much more diverse than subsemigroups.  Nevertheless, it is convenient to think about most theoretical aspects (such as extension to $\beta S$) with a full operation in mind.

The semigroups considered in this paper are $( \matn , + ) $ (with $\matn := \omega \setminus \{ 0 \}$), $(\matz,+)$ and the most important adequate partial semigroup $\matf$.

\begin{Definition}\label{def:matf}\Marginnote{ \normalsize The partial semigroup $\mathbb{F}$}
On \index{semigroup!$\matf$} \index{$\matf$} $\matf := \{ s\subseteq \omega\ |\ \nil \neq s \mbox{ finite} \}$
we define a partial semigroup structure by
\[
s \mal t := s \cup t \mbox{ if and only if } s \cap t = \nil.
\]
\end{Definition}
The theory of the Stone-\v Cech compactification allows for the (somewhat unique) extension of any operation on $S$ to its compactification, in particular a semigroup operation.
\begin{Definition}\label{def:multiplication in beta S}\index{$\beta S$!the semigroup}\Marginnote{\normalsize The semigroup $\beta S$}
For a semigroup $(S,\cdot)$, $s\in S$ and $A \subseteq S$, $p,q \in \beta S$ we define the following.\begin{itemize}
\item $ s^{-1}A := \{ t\in S\ |\ st \in A \} $.
\item $ A^{-q} := \{ s\in S\ |\ s^{-1}A\in q\} $.
\item  $ p \cdot q : = \{ A \subseteq S \ |\ A^{-q} \in p \} $. \\
Equivalently, $p\cdot q$ is generated by sets $\bigcup_{v\in V} v\cdot W_v$ for $V\in p$ and each $W_v  \in q$.
\item $ A^\star := A^{-q} \cap A$. \\
This notation will only be used when there is no confusion regarding the chosen ultrafilter.
\end{itemize}
\end{Definition}
As is well known, this multiplication on $\beta S$ is well defined and extends the operation on $S$. It is associative and right topological, i.e., the operation with fixed right hand side is continuous. For these and all other theoretical background we refer to \cite{HindmanStrauss}.

In the case of a partial semigroup, ultrafilters in $\delta S$ in a way multiply as if the partial operation was total. With the arguments from the following proposition it is a simple but useful exercise to check that if $(S, \cdot)$ is partial the above definitions still work just as well in the sense that $s^{-1}A := \{ t \in \sigma(s) \ | \ st \in A\}$ and $p\cdot q$ is only defined if it is an ultrafilter.
 
\begin{Proposition}\label{prop:semigroup delta S}\Marginnote{\normalsize The semigroup $\delta S$}
Let $S$ be a partial subsemigroup of a semigroup $T$. Then $\delta S$ is a subsemigroup of $\beta T$.
\end{Proposition}
\begin{proof} \plist
\item Simply observe that for $a \in S$
\[
\bigcup_{b \in \sigma(a)} b \cdot ( \sigma ( ab ) \cap \sigma(b) ) \subseteq \sigma(a).
\]
\item Therefore $ \sigma(S) \subseteq p \cdot q$ whenever $p,q \in \delta S$.
\pliste \end{proof} 
It is easy to similarly check that partial semigroup homomorphisms extend to full semigroup homomorphisms on $\delta S$.

Since $A^{-q}$ is not an established notation, the following useful observations present a good opportunity to test it.

\begin{Proposition}\Marginnote{\normalsize Tricks with $A^{-q}$}
Let $p,q \in \beta S$, $A\subseteq S$ and $s, t\in S$.
\begin{itemize}
\item $t^{-1}s^{-1} A= (st)^{-1} A$.
\item $s^{-1}A^{-q} = (s^{-1}A)^{-q}$.
\item $(A \cap B)^{-q} = A^{-q} \cap B^{-q}$.
\item $( s^{-1} A )^\star = s^{-1} A^\star$ (with respect to the same ultrafilter).
\item $( A^{ - q } )^{ - p } = A^{ - ( p \cdot q ) }$.
\end{itemize}
\end{Proposition}
\begin{proof} \plist
This is straightforward to check.
\pliste \end{proof} 

The proverbial big bang for the theory of ultrafilters on semigroups is the following theorem.

\begin{Theorem}[Ellis-Numakura Lemma] \label{thm:ellis-numakura} \index{Ellis-Numakura Lemma}\index{Lemma!Ellis-Numakura}
If $(S,\cdot)$ is a compact, right topological semigroup then there exists an \emph{idempotent} element in $S$, i.e., an element $p \in S$ such that $p\cdot p = p$.
\end{Theorem}
\begin{proof}
See, e.g., \cite[notes to Chapter 2]{HindmanStrauss}.
\end{proof} 

Therefore the following classical fact is meaningful.

\begin{Lemma}[Galvin Fixpoint Lemma]\label{lem:galvin}\index{lemma!{Galvin Fixpoint}}
For idempotent $p\in \beta S$, $A\in p$ implies $A^\star \in p$ and $(A^\star)^\star = A^\star$.
\end{Lemma}

\begin{proof}
$(A^\star)^\star = A^\star \cap (A^\star)^{-p} = A^\star \cap (A \cap A^{-p})^{-p} = A^\star \cap A^{-p} \cap A^{- p\cdot p} = A^\star \cap A^{-p} = A^\star$. 
\end{proof} 

The following definitions are central in what follows. Even though we mostly work in $\mathbb{N}$ and $\mathbb{F}$ we formulate them for a general setting.

\begin{Definition}\label{def:FP-sets etc} \Marginnote{\normalsize $FP$-sets, $\mathbf{x}$-support and condensations} \index{$FP$-set}
Let $\mathbf{x} = (x_n)_{n < N}$ (with $N \leq \omega $) be a sequence in a partial semigroup $(S,\mal)$ and let $K\leq \omega$.
\begin{itemize} 
\item The set of finite products (the \emph{$FP$-set}) is defined as 
\[
FP( \mathbf{x} ) := \{ \prod_{i\in v} x_i\ |\ v\in \matf \},
\]
where products are in increasing order of the indices. In this case, all products are assumed to be defined.\footnote{Note that we will mostly deal with commutative semigroups so the order of indices is not too important in what follows.}

\item $\mathbf{x}$ has \emph{unique representations} if for $v,w \in \matf$ the fact $\prod_{i \in v} x_i = \prod_{j \in w} x_j$ implies $v = w$.

\item If $\mathbf{x}$ has unique representations and $z \in FP(\mathbf{x})$ we can define the \emph{$\mathbf{x}$-support of $z$}, short $\mathbf{x}\xsupp(z)$, \index{$\mathbf{x}\xsupp$} by the equation $z = \prod_{j \in \mathbf{x}\xsupp(z)} x_j$. We can then also define $\mathbf{x}\xmin := \min \circ \mathbf{x}\xsupp$, $\mathbf{x}\xmax := \max \circ \mathbf{x}\xsupp$.


\item A sequence $\mathbf{y}=(y_j)_{j<K}$ is called a \emph{condensation}\index{condensation} of $\mathbf{x }$, in short $\mathbf{y} \sqsubseteq  \mathbf{ x } $, if
\[
FP(\mathbf{y}) \subseteq  FP( \mathbf{x} ).
\]
In particular, $\{ y_i \ |\ i<K \} \subseteq FP(\mathbf{x})$. For convenience, $\mathbf{x}\xsupp(\mathbf{y}) := \mathbf{x}\xsupp [ \{ y_i \ | \ i \in \omega \} ]$.

\item Define $FP_k( \mathbf{x} ) := FP( \mathbf{x'} )$ where $x_n' = x_{n+k}$ for all $n$.


\item $FP$-sets have a natural partial subsemigroup structure induced by $\matf$, i.e., $( \prod_{i\in s} x_i) \cdot (\prod_{i\in t} x_i)$ is defined as in $S$ but only if $\max(s) < \min(t)$. With respect to this restricted operation define $FP^\infty ( \mathbf{x } ) := \delta FP( \mathbf{x} ) = \bigcap_{k \in \omega} \bbar{ FP_k(\mathbf{x}) }$.

\item If the semigroup is written additively, we write $FS( \mathbf{x } )$ etc. accordingly (for finite sums); for $\matf$ we write $FU( \mathbf{x } )$ etc. (for finite unions).
\end{itemize}
\end{Definition}
Instead of saying that a sequence has certain properties it is often convenient to say that the generated $FP$-set does.

The following classical result is the starting point for most applications of algebra in the Stone-\v Cech compactification. We formulate it for partial semigroups.

\begin{Theorem}[Galvin-Glazer Theorem]\label{thm:galvin-glazer} \index{Galvin-Glazer Theorem}\index{Theorem!Galvin-Glazer}
Let $(S,\cdot)$ be a partial semigroup, $p\in \delta S$ idempotent and $A \in p$.
Then there exists $ \mathbf{ x } =  ( x_i )_{i\in \omega}$ in $A$ such that
\[
FP( \mathbf{x } ) \subseteq A.
\]
\end{Theorem}
\begin{proof} 
This can be proved essentially just like the the original theorem, cf.~\cite[Theorem 5.8]{HindmanStrauss}, using the fact that $\sigma(S) \subseteq p $ to guarantee all products are defined.
 \end{proof} 

An immediate corollary is, of course, the following classical theorem, originally proved combinatorially for $\matn$ in \cite{Hindman74}.

\begin{Theorem}[Hindman's Theorem] \label{thm:hindman} \index{Hindman's Theorem}\index{Theorem!Hindman} \index{Finite Sums Theorem|see{Hindman's Theorem}} \index{Theorem!Finite Sums|see{Theorem!Hindman}}
Let $S = A_0 \cup A_1$. Then there exists $i\in \{0,1\}$ and a sequence $ \mathbf{x} $ such that $FP( \mathbf{x } ) \subseteq A_i$.
\end{Theorem}

\section{Union and Strongly Summable Ultrafilters}\label{ch:union ufs}

The first part of this paper deals primarily with ultrafilters on the partial semigroup $\matf$. The following three kinds of ultrafilters were first described in \cite{Blass87-1}.

\begin{Definition}[Ordered, stable, union ultrafilters]\Marginnote{\normalsize Union Ultrafilters}
An ultrafilter $u$ on $\matf$ is called 
\begin{itemize}\index{ultrafilter!union}\index{ultrafilter!{ordered union}}\index{ultrafilter!{stable union}}
\item \emph{union} if it has a base of $FU$-sets (from disjoint sequences).
\item \emph{ordered union} if it has a base of $FU$-sets from ordered sequences, i.e., sequences $\mathbf{s}$ such that $\max(s_i) < \min(s_{i+1})$ (for all $i \in \omega$).
\item \emph{stable union} if it is union and whenever $\matf^2_< := \{ (v,w) \in \matf^2 \ | \ \max(v) < \min(w) \}$ is partitioned into finitely many pieces, there exists homogeneous $A\in u$, i.e., $A^2_<$ is included in one part.
\end{itemize}
\end{Definition}
The original definition of stability is similar to that of a P-point (or $\delta$-stable ultrafilter) which we discuss later. For their equivalence see \cite[Theorem 4.2]{Blass87-1} and \cite[Theorem 4.13]{Krautzberger09}

It is clear yet important to note that $FU$-sets always have unique representations and that all products are defined. At this point it might be useful to check the following. Union ultrafilters are elements of $\delta \matf$ and they are idempotent since for each included $FU$-set they contain all $FU_k$-sets. It is also worth while to check that if our operation on $\mathbb{F}$ was not restricted to disjoint but ordered unions then $\sigma(\mathbb{F})$ and hence $\delta \mathbb{F}$ would remain the same.

The following notion was introduced in \cite{BlassHindman87} to help differentiate union ultrafilters; it is a special case of isomorphism, but arguably the natural notion for union ultrafilters.

\begin{Definition}[Additive isomorphism]\label{def:add. isomorphism}
Given partial semigroups $S,T$, call two ultrafilters $p \in \beta S, q\in \beta T$ \emph{additively isomorphic} \index{additive isomorphism} if there exist $FP( \mathbf{ x }) \in p, FP(\mathbf{ y} )\in q$ both with unique products such that the following map maps $p$ to $q$
\[
\phii : FP( \mathbf{ x }) \rightarrow  FP(\mathbf{ y }), \prod_{i \in s} x_i \mapsto \prod_{i  \in s } y_i .
\]
We call such a map a \emph{natural (partial semigroup) isomorphism}. It extends to a homomorphism (in fact, isomorphism) between $FP^\infty( \mathbf{x} ) $ and $FP^\infty( \mathbf{y} ) $.
\end{Definition}
In the semigroup $(\matn, +)$, our interest lies in strongly summable ultrafilters.

\begin{Definition}[Strongly summable ultrafilters]
An ultrafilter $p$ on $\matn$ is called
\emph{(strongly) summable} if it has a base of $FS$-sets.
\end{Definition}
The following properties are well known and necessary to switch between summable and union ultrafilters; they are the basic tools for handling strongly summable ultrafilters, cf.~\cite{BlassHindman87}, \cite[Chapter 12]{HindmanStrauss}.
\begin{Proposition}[and Definition]\label{prop:growth condition}
Every strongly summable ultrafilter has a base of $FS(\mathbf{x})$-sets with the property \Marginnote{\normalsize Sufficient growth}
\[
(\forall n<\omega)\ x_n > 4\cdot \sum_{i<n} x_i.
\]
In this case $\mathbf{x}$ is said to have \emph{sufficient growth} which implies the following:
\begin{itemize}
\item $\sum_{i\in s} x_i = \sum_{i\in t} x_i$ iff $s= t$ (unique represenations)

\item  $\sum_{i\in s} x_i + \sum_{i\in t} x_i \in FS(\mathbf{x})$ iff $s \cap t = \nil $ (unique sums)\label{unique sums}

In particular, condensations of $\mathbf{x}$ have pairwise disjoint $\mathbf{x}$-support and the map $\sum_{i\in s} x_i \mapsto s$ maps the strongly summable to a union ultrafilter.


\item To have sufficient growth is \emph{hereditary for condensations}, i.e., if $\mathbf{x}$ has sufficient growth, so does $\mathbf{y} \sqsubseteq \mathbf{x}$ (assuming that $\mathbf{y}$ is increasing).
\end{itemize}
\end{Proposition}
\begin{proof}
This follows (in order) from \cite[Lemma 12.20, Lemma 12.34, Lemma 12.32, Theorem 12.36]{HindmanStrauss}. The last observation follows easily from the second bullet and the growth of $\mathbf{x}$ since the growth of $\mathbf{x}$ implies that to be increasing means to be $\mathbf{x}\xmax$-increasing.
\end{proof} 

Maybe the most important aspect to remember is this: whenever we have a condensation of a sequence with sufficient growth, its elements have pairwise disjoint $\mathbf{x}$-support (by \ref{unique sums}.2) and we can apply the much less messy intuition about $FU$-sets to understand the structure of the $FS$-set. In particular, whenever a sequence $\mathbf{x}$ in $\matn$ has sufficient growth we can apply the terminology of $\mathbf{x}\xsupp$, $\mathbf{x}\xmax$ and $\mathbf{x}\xmin $ as introduced in the preliminaries.

Although it is not relevant in our setting note that on the one hand growth by a factor $2$ (instead of $4$) already implies the above properties (with identical proofs as in the references). On the other hand the proof of \cite[Lemma 12.20]{HindmanStrauss} can easily be enhanced to show that for any $k \in \mathbb{N}$ every strongly summable ultrafilter will have a base with growth factor $k$ which leads to other interesting properties such as \cite[Lemma 12.40]{HindmanStrauss}.

\section{Strongly summable ultrafilters are special}

Recall that we aim to extend a theorem by N.~Hindman and D.~Strauss on writing strongly summable ultrafilters as sums originally published in \cite{HindmanStrauss95}, cf.~\cite[Theorem 12.45]{HindmanStrauss}. The original result was shown for a certain class of strongly summable ultrafilters, the so-called special strongly summable ultrafilters. Our main result will extend this to a wider class of strongly summable ultrafilters. The proof will require one new observation, which we prove in this section, as well as a series of modifications of the original proof as presented in, e.g., \cite[Chapter 12]{HindmanStrauss}.

To investigate special strongly summable ultrafilters as described in \cite{HindmanStrauss95} and \cite[12.24]{HindmanStrauss}, it is useful to switch to union ultrafilters. However, the notion introduced below is strictly weaker than the original one used by N.~Hindman and D.~Strauss.

\begin{Definition}\label{def:special summable}
Let $\mathbf{x}, \mathbf{y}$ be sequences in $\matn$.
\begin{itemize}
\item A strongly summable ultrafilter \Marginnote{Special strongly summable ultrafilter} $p \in \beta \matn$ is \emph{special} \index{ultrafilter!{special summable}} if there exists $FS(\mathbf{x}) \in p$ with sufficient growth such that
\[
(\forall L \in [\omega]^\omega) ( \exists \mathbf{y} \sqsubseteq \mathbf{x})\ FS(\mathbf{y}) \in p \mbox{ and } | L \setminus \mathbf{x}\xsupp(\mathbf{y}) | = \omega.
\]
Given the sequence $\mathbf{x}$ we say that $p$ is special \emph{with respect to $\mathbf{x}$}.

\item A union ultrafilter \Marginnote{Special union ultrafilter}\index{ultrafilter!{special union}} $u\in \beta \matf$ is \emph{special} if
\[
( \forall L \in [\omega]^\omega ) (\exists X \in u ) |L \setminus \bigcup X | = \omega.
\]
\end{itemize}
\end{Definition}

In \cite{HindmanStrauss95} and \cite[Chapter 12]{HindmanStrauss}, the notion of ``special'' is in this terminology ``special with respect to $(n!)_{n\in \omega}$ and additionally \emph{divisible}'', i.e., there is a base of sets $FS(\mathbf{x})$ with $x_n | x_{n+1}$ for all $n\in \omega$. However, \cite[Theorem 5.8]{HindmanStrauss95} gives an example of a strongly summable ultrafilter that is not additively isomorphic to a divisible ultrafilter so our notion is consistently weaker.

It is not surprising yet very useful that to be the witness for specialness is hereditary for condensations.

\begin{Proposition}\label{prop:special hereditary}\Marginnote{\normalsize Special is hereditary}
If a strongly summable ultrafilter $p$ is special with respect to $\mathbf{x}$ and $\mathbf{y} \sqsubseteq \mathbf{x}$ with $FS(\mathbf{y}) \in p$, then $p$ is special with respect to $\mathbf{y}$.
\end{Proposition}
\begin{spoiler}
The uniqueness of $\mathbf{x}$-support allows us to link the elements of the $\mathbf{y}$-support to the $\mathbf{x}$-support. Hence, for a common condensation, missing elements in the $\mathbf{x}$-support will imply missing elements in the $\mathbf{y}$-support.
\end{spoiler}

\begin{proof} \plist
\item Take any $L \in [ \omega ]^{ \omega  }$.

\item Then define $L': = \{ i \in \omega \ | \ (\exists k \in L )\ i \in \mathbf{x} \xsupp (y_k)  \} $. $L'$ is obviously infinite.
\begin{prose}
Note that the $k$'s are unique thus linking the two kinds of support.
\end{prose}

\item Since $p$ is special there exists a condensation $\mathbf{z} \sqsubseteq \mathbf{x}$ with $FS(\mathbf{z}) \in p$ and
\[
| L' \setminus \mathbf{x}\xsupp(\mathbf{z}) | = \omega.
\]

\item For a common condensation $\mathbf{v} \sqsubseteq \mathbf{y}, \mathbf{z}$ with $FS(\mathbf{v}) \in p$, naturally
\[
| L' \setminus \mathbf{x}\xsupp(\mathbf{v}) | = \omega.
\]

\item $| L \setminus \mathbf{y}\xsupp(\mathbf{v}) | = \omega$.
\plist
\item If $i \in L' \setminus \mathbf{x}\xsupp(\mathbf{v})$, then there exists (by definition of $L'$) some $k_i \in L$ with $i \in \mathbf{x} \xsupp(y_{k_i})$.

\item  But then no $v_j$ can have $k_i \in \mathbf{y}\xsupp (v_j)$ (or else $x_i \in \mathbf{x}\xsupp (v_j)$ which is impossible due to the previous proposition).

\item  In other words, $k_i \in L \setminus \mathbf{y}\xsupp(\mathbf{v})  $.

\item  Since $| L' \setminus \mathbf{x}\xsupp(\mathbf{v}) | = \omega$ and the map $i \mapsto k_i$ is finite-to-one, $| L \setminus \mathbf{y}\xsupp(\mathbf{v}) | = \omega$.
\pliste
\item This completes the proof.
\pliste \end{proof}

The second observation is that the notions of special summable and special union ultrafilters are in fact equivalent.

\begin{Proposition}\index{ultrafilter!{special union}}\index{ultrafilter!{special summable}} \Marginnote{\normalsize Special union = special strongly summable}
Let $p$ be a strongly summable ultrafilter additively isomorphic to a union ultrafilter $u$. Then $p$ is special if and only if $u$ is.
\end{Proposition}

\begin{proof} \plist
\item Assume that $p$ and $u$ are as above and additively isomorphic via
\[
\phii: FS(\mathbf{x}) \rightarrow FU(\mathbf{s}), \sum_{i\in F} x_i \mapsto F,
\]
for suitable sequences $\mathbf{x}, \mathbf{s}$ in $\mathbb{N}$ and $\mathbb{F}$ respectively.

\item By switching to a condensation we may assume that $\mathbf{x}$ has sufficient growth.

\item If $u$ is special, then $\phii$ clearly guarantees that $\mathbf{x}$ is a witness for $p$ being special.

\item If $p$ is special, we can assume that $\mathbf{x}$ is a witness of specialness thanks to the preceding proposition.
\item Then again $\phii$ will guarantee that $u$ is special.
\pliste \end{proof} 

The key fact is that all union ultrafilters are special.

\begin{Theorem}[Union ultrafilters are special] \label{thm:union ufs are special} \index{ultrafilter!{special union}}\index{ultrafilter!{special summable}}\Marginnote{\normalsize Union ultrafilters are special}
Every union ultrafilter is special. Accordingly, all strongly summable ultrafilters are special.
\end{Theorem}

\begin{spoiler}
Assuming that some set covers all of $L$, a parity argument on pairs of the form $(i,i+1)$ in the support will yield a condensation that misses a lot of $L$.
\end{spoiler}

\begin{proof} \plist
\item Let $L\in [\omega]^\omega$.
\begin{prose}
Remember that $ \bigcup FU(\mathbf{s}) = \bigcup \{ s_i \ | \ i \in \omega \}$. For the mental picture of the arguments it is helpful (though not necessary) to enumerate sequences according to the maximum.
\end{prose}
\item We may assume that $\{ s\in \mathbb{F} \ |\ s\cap L \neq \emptyset \} \in u$.
\plist
\item Otherwise its complement, call it $X$, has $L \setminus \bigcup X = L$ infinite -- as desired.
\pliste
\item  Since $u$ is a union ultrafilter, we find \hl{$FU(\mathbf{s}) \in u$} included in this set.
\item If $L \setminus \bigcup FU(\mathbf{s})$ is infinite, we are done.
\item So assume it is finite; without loss it is empty.
\begin{prose}
In the following sense we can now think as if $L= \omega$. If $\mathbf{t} \sqsubseteq \mathbf{s}$ and $i \not\in \mathbf{s}\xsupp(\mathbf{t})$, then $s_i \cap L \neq \emptyset$ but $s_i \cap \bigcup FU(\mathbf{t}) = \emptyset$. So dropping elements in the $\mathbf{s}$-support means dropping elements in $L$ (and vice versa). So we can concentrate on $\mathbf{s}\xsupp (\mathbf{s}) = \omega$.
\end{prose}
\item Consider $\hl{ \pi : FU(\mathbf{s}) \rightarrow \omega}, t \mapsto \{ i : s_i, s_{i+1} \subseteq t \}$.  We're interested in whether $\pi(t)$ is even or odd.
\item Since $u$ is a union ultrafilter, we can find \hl{$FU(\mathbf{t}) \in u$} such that the elements of $\pi[FU(\mathbf{t})]$ all have the same parity.
\item \Marginnote{The parity argument.}But the elements of $\pi[FU(\mathbf{t})]$ can only be of even size.
\plist
\item For any $x\in FU(\mathbf{t})$, there exists $i,j \in \omega$ such that $\mathbf{s}\xmax (x) < i < \mathbf{s}\xmin(t_j)$.
\item In that case $\pi( x \cup t_j) = \pi(x) + \pi(t_j)$ -- which is even since $\pi(x) = \pi(t_j)$.
\pliste
\item Then $L \setminus \bigcup FU(\mathbf{t})$ is infinite.
\plist
\item Assume towards a contradiction that it is finite. 
\begin{prose}
We will study the gaps in the $\mathbf{s}$-support of elements in $FU(\mathbf{t})$ since they correspond to elements in $L \setminus FU(\mathbf{t})$.
\end{prose}
\item The set $\mathbf{s}\xsupp(\mathbf{t})$ must be cofinite since $\mathbf{s}$ covers all of $L$ and every $s_i \cap L \neq \emptyset$.
\item In other words, there exists \hl{$b\in \omega$} such that $(\forall i \geq b) (\exists j_i) s_i \subseteq t_{j_i}$.
\begin{prose}
Consider for a moment $t_{j_b}$, the $t_j$ containing $s_b$. Since $\mathbf{t}$ covers all later $s_i$, some $t_j$ contains $\mathbf{s}\xmax(t_{j_b}) +1 $. Therefore their union ``gains'' a pair of adjacent indices, i.e.,  $\pi(t_{j_b} \cup t_j) \geq \pi(t_{j_b}) + \pi(t_j) +1$. Since $\pi(t_{j_b} \cup t_j)$ is even it must ``gain'' even more. If $\mathbf{t}$ was ordered, this would be impossible. For the unordered case, we need to argue more subtly.
\end{prose}
\item We define \hl{$x := \bigcup_{ i \leq b} t_{j_i} \in FU(\mathbf{t})$}, adding to $t_{j_b}$ everything ``below'' it.
\begin{prose}
$x$ is our initial piece. It contains the $\mathbf{s}\xsupp(\mathbf{t})$ up to $b$. This ensures that any $t_j$ disjoint from $x$ must have $\mathbf{s}$-support beyond $b$.
\end{prose}
\item Next we define \hl{ $b_1 := \mathbf{s}\xmax(x)$}, i.e., the index of the last $s_i \subseteq x$.
\item Of course,  $b_1 \geq b$ by choice of $t_{j_b} \subseteq x$.
\begin{prose}
We will derive the contradiction from the fact that we can fill the entire interval $[b, b_1+1]$ by choice of $b$.
\end{prose}
\item Then we define \hl{$y := t_{j_{b_1+1}}$}, i.e., the $t_j$ that contains the next element of the $\mathbf{s}$-support.
\item Finally, let \hl{$z :=( \bigcup \{ t_{j_i} : i <b_1\}) \setminus  (x\cup y) $}.
\begin{prose} $y$ follows on where $x$ ends, $z$ fills all the gaps in the $\mathbf{s}$-support of $x \cup y$ between $b$ and $b_1$ (and, of course, the support of $z$ lies only beyond $b$). We will now analyze how gaps in $\mathbf{s}\xsupp(x)$ are actually filled.
\end{prose}
\item On the one hand, we can compare $\pi(x)$ and $\pi(x \cup y)$.
\item By definition,
\[
\pi(x\cup y) = \pi(x) \dot\cup \pi(y) \dot\cup \{ i \ |\ s_i \subseteq x, s_{i+1} \subseteq y \mbox{ or vice versa} \}.
\]
Let us call elements in the third set \emph{emerged} indices.
\item We know that $\pi( x \cup y)$ contains one emerged index, namely $b_1$.
\item But $\pi(x \cup y)$ is even and $x$ has no support past $b_1$. 
\item Therefore $\pi(x \cup y)$ must have an odd number of emerged indices below $b_1$.
\item In particular, $y$ has $\mathbf{s}$-support below $b_1$ (sitting inside the gaps of the $\mathbf{s}$-support of $x$).
\item\Marginnote{Four types of gaps}There are four ways how those $i \in \mathbf{s}\xsupp (y)$ with $i< b_1$ can be found within the gaps of $\mathbf{s}\xsupp(x)$: only at the beginning of a gap, only at the end of a gap, both at the beginning and end of a gap and finally at neither beginning nor end of a gap.
\item The latter two cases do not change the parity of $\pi(x\cup y)$ since they account for two and zero emerged indices respectively.
\item So to make up for $b_1$ there must be an odd number of cases where $\mathbf{s}\xsupp(y)$ fills only the beginning or only the end of a gap in $\mathbf{s}\xsupp (x)$.
\item On the other hand, we can similarly compare $\pi(x)$ and $\pi(x \cup z)$.
\item We know that $\mathbf{s} \xsupp( x\cup y \cup z)$ contains the entire interval $[b,b_1]$.
\item In particular, $\mathbf{s}\xsupp(z)$ fills the beginning or end of any gap of $\mathbf{s}\xsupp(x)$ that was not filled by $\mathbf{s}\xsupp(y)$.
\item By the above analysis of $\pi(x\cup y)$ and $\pi(x)$ this gives an odd number of emerged indices in $\pi(x \cup z)$ below $b_1$.
\item But then $\pi(x \cup z)$ is odd since $z$ has no support below $b$, $x$ has no support above $b_1$ and neither contains $b_1+1$ -- a contradiction.
\pliste
\pliste \end{proof} 

I am very grateful for Andreas Blass's help in closing a gap in the final step of the above proof.

\section{Disjoint support and trivial sums}
  
There is need for another notion of support before formulating the main result. Every divisible sequence $\mathbf{a} = (a_n)_{n \in \omega }$, i.e., with $a_n | a_{n+1}$ for $n \in \omega$, with $a_0 = 1$ induces a unique representation of the natural numbers; the easiest case to keep in mind would be $a_n = 2^n$, i.e., the binary representation. We will work with an arbitrary divisible sequence but it might be best to always think of the binary case.

\begin{Definition}\Marginnote{\normalsize Fix $\mathbf{a}$, the divisible sequence} 
For the rest of this section we fix some divisible sequence $\mathbf{a} = (a_n)_{n \in \omega }$, i.e., with $a_n | a_{n+1}$ for $n \in \omega$, with $a_0 = 1$.
\begin{itemize}
\item We consider $\prod_{i \in \omega} \frac{a_{i+1}}{a_i} = \prod_{i \in \omega} \{0,\ldots, \frac{a_{i+1}}{a_i} -1 \}$ as a compact, Hausdorff space (with the product topology, each coordinate discrete).

\item We can then define $\alpha:\matn \rightarrow \prod_{i \in \omega} \frac{a_{i+1}}{a_i}$ by the (unique) relation 
\[
n = \sum_{i\in \omega} \alpha(n)(i) \cdot a_i.
\]
In other words, $\alpha(n)$ yields the unique representation of $n$ with respect to $\mathbf{a}$. Note that $\alpha(n)$ has only finitely many non-zero entries for any $n$ but for $p \in \beta \mathbb{N}$ its continuation $\alpha(p)$ might not.

\item The \emph{$\alpha$-support of $n$}, $\alpha\xsupp(n)$, is the (finite) set of indices $i$ with $\alpha(n)(i) \neq 0$; similarly we define $\alpha\xmax(n)$, $\alpha\xmin$ to be its maximum and minimum respectively.

\item A sequence $\mathbf{x}=(x_n)_{n \in \omega}$ has \emph{disjoint $\alpha$-support} if its elements do; allowing confusion, $FS(\mathbf{x})$ is said to have disjoint support.

\item A strongly summable ultrafilter has \emph{disjoint $\alpha$-support} if it contains an $FS$-set with disjoint $\alpha$-support and sufficient growth.

\item An idempotent ultrafilter $p$ can be \emph{written as a sum only trivially} if
\[
(\forall q,r\in \beta \matn)\ q+r = p \Rightarrow q,r \in (\matz+p)
\]

\item For $(2^n)_{n\in \omega}$, the binary support is abbreviated $\bsupp$; its maximum and minimum by $\bmax$ and $\bmin$ respectively.
\end{itemize}
\end{Definition}

For the ``trivial sums'' property we should note that it is an easy exercise to show that $\beta \matn \setminus \mathbb{N}$ is a left ideal of $(\beta \matz,+)$; in particular $\mathbb{Z} +p \subseteq \beta \mathbb{N}$.

So far we have always been interested in the finite sums of a sequence. It might therefore cause confusion as to why we chose the $\alpha$-support when we have so far only studied the $\mathbf{a}$-support (which only coincides on $FS(\mathbf{a})$). Why not just assume that $FS(\mathbf{a})$ is in our strongly summable ultrafilter? From a certain point of view, this is what happens in the original result by Hindman and Strauss, cf.~\cite[12.24]{HindmanStrauss} and in \cite{HindmanStrauss95}. The advantage of our notion of disjoint $\alpha$-support lies precisely in dropping this requirement -- we won't need (a suitable condensations of) $FS(\mathbf{a})$ in the strongly summable ultrafilter. In this spirit, there hopefully won't be a lot of confusion between $\alpha$-support and $\mathbf{a}$-support. Nevertheless we will see that the reasoning with $\alpha$-support is quite similar when considering sequences with disjoint $\alpha$-support.

Since we will be concerned with $\bigcap_{n\in \matn} \bbar{a_n \matn}$ it is worthwhile to point out that by divisibility, $a_n\matn \supseteq a_{n+1}\matn$. Therefore an ultrafilter containing infinitely many such sets already contains all of them. Also, it is well known that any idempotent ultrafilter contains the set of multiples for any number. The following will be the main result.

\begin{Theorem}[Strongly summable ultrafilters as sums]\Marginnote{\normalsize Trivial Sums}
Every  strongly summable ultrafilter with disjoint $\alpha$-support can be written as a sum only trivially.
\end{Theorem}

The proof requires a series of technical propositions, but the following convenient corollary is immediate.

\begin{Corollary}
Every strongly summable ultrafilter is additively isomorphic to a strongly summable ultrafilter that can only be written as a sum trivially.
\end{Corollary}
\begin{proof} \plist
\item For any strongly summable ultrafilter $p$, pick $FS(\mathbf{x}) \in p$ with sufficient growth. 

\item Then, e.g., the natural additive isomorphism $\phii$ between $FS(\mathbf{x})$ and $FS((2^n)_{n \in \omega})$ maps $p$ to a strongly summable ultrafilter with disjoint binary support.
\plist
\item Let $p'$ be the image of $p$; clearly, $p'$ is a strongly summable ultrafilter.
\item Fix some $FS(\mathbf{y}) \in p'$ with sufficient growth. 
\item Then $FS(\mathbf{y}) = \phii[FS(\mathbf{z})] = FS(\phii[\mathbf{z}])$ for some $\mathbf{z} \sqsubseteq \mathbf{x}$.
\item The growth of $\mathbf{x}$ guarantees that each $z_i$ is a disjoint union of elements from $\mathbf{x}$.
\item Hence each $y_i$ is a disjoint union of $\phii[\mathbf{x}] = (2^n)_{n\in \omega}$.
\item In other words, $\mathbf{y}$ has disjoint binary support, as desired.
\pliste 
\pliste \end{proof} 

In \cite{HindmanStrauss95} it is shown that strongly summable ultrafilters that are divisible and special with respect to $(n!)_{n\in \omega}$ can only be written as a sum trivially; however, by \cite[Theorem 5.8]{HindmanStrauss95}, there consistently exist strongly summable ultrafilters that are not additively isomorphic to a divisible strongly summable ultrafilter.\footnote{cf.~the comment after Definition \ref{def:special summable}. We could summarize our approach as replacing $(n!)_{n\in \omega}$ with $\mathbf{a}$ and divisibility with disjoint $\alpha$-support.} In so far, this is an improvement.

To begin the series of technical observations, note one additional detail concerning the herditary nature of specialness.

\begin{Lemma}\label{lem:summable disjoint binary support and special}\index{ultrafilter!{special summable}}\Marginnote{\normalsize $\alpha$-special}
A strongly summable ultrafilter $p$ with disjoint $\alpha$-support is also $\alpha$-special in the sense that there exists $FS(\mathbf{x}) \in p$
\[
(\forall L \in [\omega]^\omega) ( \exists \mathbf{y} \sqsubseteq \mathbf{x})\ FS(\mathbf{y}) \in p \mbox{ and } | L \setminus \alpha\xsupp(\mathbf{y}) | = \omega.
\]
\end{Lemma}
\begin{spoiler}
We argue as for the heredity of specialness using a common condensation of witnesses for disjoint $\alpha$-support and specialness.
\end{spoiler}

\begin{proof} \plist
\item Pick $\mathbf{x}$ as a witness for the disjoint $\alpha$-support of a strongly summable ultrafilter $p$.

\item We may assume that $\mathbf{x}$ also witnesses that $p$ is special.
\plist
\item By Proposition \ref{prop:special hereditary}, to be the witness for specialness is hereditary.
\item By Proposition \ref{prop:growth condition}, any condensation of $\mathbf{x}$ has pairwise disjoint support $\mathbf{x}$-support, hence pairwise disjoint $\alpha$-support; in other words, to have disjoint $\alpha$-support is hereditary.
\item Therefore a common condensation of the respective witnesses will have both properties.
\pliste
\item Given $L \in [\omega]^\omega$; if $L \setminus \alpha\xsupp(\mathbf{x})$ is infinite, we are done.

\item If not we can consider the (infinite) set
\[
L':= \{ n\ | \ (\exists i \in L)\ i \in \alpha\xsupp(x_n)\}.
\]
\item By specialness there exists $\mathbf{y} \sqsubseteq \mathbf{x}$ with $FS(\mathbf{y}) \in p$ and $L' \setminus \mathbf{x}\xsupp(\mathbf{y})$ infinite.

\item But this implies $L\setminus \alpha\xsupp(\mathbf{y})$ is infinite by choice of $L'$ and the disjoint $\mathbf{x} \xsupp$ of members of $\mathbf{y} $.
\pliste \end{proof}

The following well known theorem proves, in a manner of speaking, half the theorem.
\begin{Theorem}\label{cor:summables strongly right maximal}\index{strongly right maximal}
Every strongly summable ultrafilter $p$ is a \index{strongly right maximal}\emph{strongly right maximal idempotent}, i.e., the equation $q + p = p$ has the unique solution $q = p$.
\end{Theorem}
\begin{proof}
This is, e.g., \cite[Theorem 12.39]{HindmanStrauss}.
\end{proof} 
The next result is also well known and easily checked.

\begin{Proposition} \label{prop: nearly prime semigroups}
For $n \in \matn$, $q,r\in \beta \matn$ the following holds.
\begin{itemize}
\item If $q+r \in \bbar{ n \matn } $, then either both $q, r \in \bbar{ n \matn} $ or neither is.
\item 
Similarly we can replace $\bbar{n\matn}$ by $\bigcap_{n\in \matn} \bbar{a_n \matn}$ and $\matz + \bigcap_{n\in \matn} \bbar{a_n \matn}$.
\end{itemize}
\end{Proposition}
\begin{proof}
This is, e.g., \cite[Lemma 2.6]{HindmanStrauss95}.
\end{proof}

As mentioned earlier, our proof follows the same strategy as the proof in \cite{HindmanStrauss95} and \cite[Chapter 12]{HindmanStrauss}; the proof for the right summand consists of two parts. The first part proves that if one of the summands is close to the strongly summable ultrafilter, i.e., in $\bigcap_{n \in \omega} \bbar{a_n\matn}$, it is already equal. The second part shows that writing a strongly summable ultrafilter with disjoint support as a sum can only be done with the summands ``close enough'' to it.

For the first part, a technical lemma reflects the desired property: under restrictions typical for ultrafilter arguments, elements of an $FS$-set with disjoint $\alpha$-support can be written as sums only trivially.

\begin{Lemma}[Trivial sums for $FS$-sets] \label{lem:trivial sums for FS-sets}
Let $\mathbf{x}=(x_n)_{n\in \matn}$ be a sequence with disjoint $\alpha$-support and enumerated with increasing $\alpha\xmin$, $a\in \matn$ and
\[
m:= \min \{ i\ |\ \alpha\xmax(a) < \alpha\xmin(x_i)\},
\]
Then for every $b \in \matn$ with $\alpha\xmax(x_m) < \alpha\xmin(b)$
\[
a + b \in FS(\mathbf{x}) \Rightarrow a, b \in FS(\mathbf{x}).
\]
\end{Lemma}

\begin{spoiler}
The simple idea is that neither the sums of the $x_i$ nor the sum $a+b$ will have any carrying over in the $\alpha$-support. Hence, the $\mathbf{x}$-support of $a+b$ splits into $\mathbf{x}$-support of $a$ and $b$.
\end{spoiler}

\begin{proof} \plist
\item Assume $\mathbf{x}$, $a$ and $b$ are given as in the lemma.

\item Since $a+b \in FS(\mathbf{x})$, there exists some finite, non-empty $H \subseteq \matn$
\[
\circledast \qquad a+b = \sum_{i\in H} x_i.
\]

\item Define
\[
H_a := \{ j \in H \ |\ \alpha\xsupp(x_j) \cap \alpha\xsupp(a) \notnil\}
\]
and $H_b$ similarly.

\item $H = H_a \dotcup H_b$.
\plist
\item On the one hand $\alpha\xsupp(a) \cap \alpha\xsupp(b) = \nil$ by assumptions on $b$; also $\mathbf{x}$ has disjoint $\alpha$-support.

\item So there is no carrying over (in the $\alpha$-support) on either side of the equation $\circledast$, i.e.,
\[
H = H_a \cup H_b.
\]
\item On the other hand, if $\alpha\xsupp(x_i) \cap \alpha\xsupp(a) \notnil$, then $i\leq m$ by choice of $m$.

\item This in turn implies $\alpha\xsupp(x_i) \cap \alpha\xsupp(b) = \nil$ by the choice of $b$.
\item In other words, $H_a \cap H_b = \nil $.
\pliste
\item Then $\sum_{i\in H_a} x_i = a$ and $\sum_{i\in H_b} x_i = b$ --  as desired.
\pliste \end{proof}

The next lemma takes the proof nearly all the way, i.e., if the second summand is ``close enough'' to the strongly summable ultrafilter, both are equal to it.

\begin{Lemma}[Trivial sums for $\bigcap_{n\in \omega} \bbar{a_n \matn}$]\label{lem:trivial sums for H}
For any strongly summable ultrafilter $p$ with disjoint $\alpha$-support
\[
(\forall q\in \beta\matn)(\forall r\in \bigcap_{n\in \omega} \bbar{a_n \matn})\ q+r= p \Rightarrow q=r=p.
\]
\end{Lemma}

\begin{spoiler}
The proof is basically a reflection argument. Arguing indirectly, the addition on $\beta \matn$ reflects to elements in the sets of the ultrafilters in such a way that non-trivial sums of ultrafilters lead to non-trivial sums of an $FS$-set, contradicting Lemma \ref{lem:trivial sums for FS-sets}.
\end{spoiler}

\begin{proof} \plist
\item Since any strongly summable ultrafilter is strongly right maximal by Theorem \ref{cor:summables strongly right maximal}, it suffices to show that $r = p$. \hl{Assume to the contrary} that $r\neq p$.

\item Pick a witness for $p$, i.e., $\mathbf{x}=(x_n)_{n\in \matn}$ with sufficient growth and disjoint $\alpha$-support; without loss $FS(\mathbf{x}) \in p \setminus r$.

\item Since $q + r = p$, $FS(\mathbf{x})^{-r} \in q$; so pick \hl{$a$} such that $-a + FS(\mathbf{x}) \in r$.

\item Pick \hl{$m$} as for Lemma \ref{lem:trivial sums for FS-sets}, i.e., such that all $(x_n)_{n>m}$ have $\alpha\xmax(a) < \alpha\xmin(x_n)$ (which is possible since $\mathbf{x}$ has disjoint $\alpha$-support).

\item Define \hl{$M:= \alpha\xmax(x_m)+1$}; note that the multiples of $a_M$ have $\alpha$-support beyond the support of both $x_m$ and $a$.

\item Now
\[
(-a+FS(\mathbf{x})) \cap (\matn\setminus FS(\mathbf{x})) \cap a_M\matn \in r.
\]
So pick \hl{$b$} from this intersection.

\item Then $a+b \in FS(\mathbf{x})$. But applying Lemma \ref{lem:trivial sums for FS-sets} both $a,b \in FS(\mathbf{x})$ \gflash contradicting $b \notin FS(\mathbf{x})$.
\pliste \end{proof}

In the final and main lemma, it remains to show that if a strongly summable ultrafilter is written as a sum, then the summands are already ``close enough''.

\begin{Lemma}[Nearly trivial sums]\label{lem:nearly trivial sums}
For any strongly summable ultrafilter $p$ with disjoint $\alpha$ support
\[
(\forall q,r\in \beta\matn)\ q+r= p \Rightarrow q,r \in \matz + \bigcap_{n\in \omega} \bbar{a_n\matn}.
\]
\end{Lemma}

\begin{spoiler}
We follow the strategy of the proof of \cite[Theorem 12.38]{HindmanStrauss} The argument is similar to the previous lemma, i.e., if $q \notin \matz + \bigcap_{n\in \omega} \bbar{a_n\matn}$, there will always be a sum $a+b$ that cannot end up in a certain $FS$-set. For this, the image of $q$ under (the continuous extension of) $\alpha$ is analyzed. Using the fact that strongly summable ultrafilters are special, it turns out that there cannot be enough carrying over available to always end up in the $FS$-set.
\end{spoiler}

\begin{proof} \plist

\item By Proposition \ref{prop: nearly prime semigroups} it suffices to show that $q \in \matz + \bigcap_{n\in \omega} \bbar{a_n\matn}$.

\item Define the following subsets of $\omega$.
\begin{align*}
& \hl{Q_0:= }\{ i \in \omega \ | \ \alpha(q)(i) < \frac{a_{i+1}}{a_i}-1 \} \\ 
& \hl{Q_1:= } \{ i \in \omega \ |\ \alpha(q)(i) > 0 \}. 
\end{align*}
\begin{prose}
In other words, $Q_0$ counts where the $\alpha$-support does not have a maximal entry, $Q_1$ counts where it does not have a minimal entry, i.e., $Q_1$ is just the support of the function $\alpha(q)$ in the usual sense.
\end{prose}

\item If either $Q_0$ or $Q_1$ is finite, then $q \in \matz + \bigcap_{n\in \omega} \bbar{a_n\matn}$.
\plist
\item \textbf{Case 1} $Q_1$ is finite.

\item Pick $k\in \omega$ such that $\alpha(q)(n) = 0$ for $n > k$.

\item Then \Marginnote{Shift by the non-trivial part of $\alpha(q)$} show that $z := \sum_{i\leq k} \alpha(q)(i) a_i$ has
\[
(\forall n>k)\ z + a_n \matn \in q,
\]
\plist
\item Given $n>k$ define
\[
U_{z,n} := \{ s\in \prod_{i\in \omega} \frac{a_{i+1}}{a_i} |\ s\restr n = \alpha(q)\restr n = \alpha(z) \restr n\}.
\]
\item Obviously, $U_{z,n}$ is an open neighbourhood of $\alpha(q)$, hence $\alpha^{-1}[U_z] \in q$.

\item But it is easily checked that $\alpha^{-1}[U_z] = z + a_n \matn$.
\pliste
\item Since $\mathbf{a}$ was divisible, $q \in z + \bigcap_{n\in \omega} \bbar{a_n\matn}$ -- as desired.

\item \textbf{Case 2} $Q_0$ is finite.

\item Pick $k$ such that $\alpha(q)(n) = \frac{a_{n+1}}{a_n}$, i.e., maximal, for $n > k$.

\item This \Marginnote{Shift by the non-trivial part of $\alpha(q)$} time show that $z:= a_{k+1} - \sum_{i<k} \alpha(q)(i) a_i$ has
\[
(\forall n>k)\ -z + a_n \matn \in q,
\]
and therefore again $q \in -z + \bigcap_{n\in \omega} \bbar{a_n\matn}$.
\plist
\item Again, given $n>k$, consider $\alpha^{-1}[U_{z,n}]$.

This time we check that $\alpha^{-1}[U_{z,n}] = -z+ a_n \matn$.

\item Let $w \in \alpha^{-1}[U_{z,n}]$. Then for some $b\geq 0$
\[
w= b \mal a_{n+1} + \sum_{i>k}^n (\frac{a_{i+1}}{a_i}-1) a_i + \sum_{i\leq k} \alpha(q)(i) a_i,
\]
since by assumption that $Q_0$ is finite, i.e., all of $\alpha(q)(i)$ beyond $k$ is maximal.

\item But this implies\Marginnote{Telescope sums}
\begin{align*}
w+z & = b \mal a_{n+1} + \sum_{i >k}^n ( \frac{a_{i+1}}{a_i} -1) a_i + a_{k+1} \\
 & = b\mal a_{n+1} + a_{n+1} = (b+1) a_{n+1},
\end{align*}
as desired.
\pliste
\item This concludes case 2.
\pliste

\item\Marginnote{Assume $Q_0, Q_1$ infinite}So let us \hl{assume to the contrary} that $q \notin \matz + \bigcap_{n\in \omega} \bbar{a_n \matn}$, i.e., both $Q_0, Q_1$ are infinite.

\item Since $u$ is strongly summable with disjoint $\alpha$-support, pick a sequence \hl{$\mathbf{x}=(x_n)_{n\in \omega}$} with disjoint $\alpha$-support, sufficient growth and $FS(\mathbf{x}) \in u$.

\item By Lemma \ref{lem:summable disjoint binary support and special}, assume without loss that both $Q_0\setminus \alpha\xsupp(\mathbf{x})$ and $Q_1\setminus \alpha\xsupp(\mathbf{x})$ are infinite.

\begin{prose}
Towards the final contradiction, it is now necessary to choose a couple of natural numbers; each choice will be followed by a short comment.
\end{prose}

\item 
By $q+r = p$ of course $FS(\mathbf{x})^{-r} \in q$; so pick \hl{$a$ with $-a+FS(\mathbf{x}) \in r$}.
\begin{prose}
$a$ can $r$-often be translated into $FS(\mathbf{x})$ -- which will be too often.
\end{prose}
\item 
Next, pick \hl{ $s_1\in Q_1\setminus \alpha\xsupp(\mathbf{x})$ and $s_2 \in Q_0 \setminus \alpha\xsupp(\mathbf{x})$ } with
\[
s_2 > s_1 > \alpha\xmax(a)
\]
\begin{prose}
On the one hand, $s_1$ ensures $\sum_{i\leq s_2} \alpha(q)(i) a_i - a > 0$, but this difference has a non-maximal entry at $\alpha\xmax$ since $s_2 \in Q_0$. On the other hand, $\alpha(q)(s_2)$ is not maximal, $\alpha(q)(s_1)$ is not minimal, but every $z\in FS(\mathbf{x})$ has $\alpha(z)(s_2) = \alpha(z)(s_1) = 0 $.
\end{prose}

\item 
By $q+r = p$ also $(a_{s_2+1}\matn)^{-r} \in q$, so \hl{pick $b$ with}
\[
b \in (a_{s_2+1}\matn)^{-r} \cap ( \sum_{i\leq s_2} \alpha(q)(i) a_i + a_{s_2+1}\matn) \in q.
\]
where the latter set is in $q$ since it is $U_{q\restr(s_2+1), s_2+1}$; cf.~Step 3.

\begin{prose}
So $b$ has $\alpha(b)(s_i) = \alpha(q)(s_i)$ (for $i=2,1$), i.e., non-maximal and non-minimal respectively.
In particular, $b-a > a_{s_1} - a > 0$ but $\alpha( b - a)(s_2)$ is not maximal.
\end{prose}

\item 
Finally, \hl{choose $y \in (-b + a_{s_2+1} \matn) \cap (-a + FS(\mathbf{x})) \in r$}.

\begin{prose}
Note that since $s_2 \notin \alpha\xsupp(\mathbf{x})$ and $a+y \in FS(\mathbf{x})$ we have $\alpha(a+y)(s_2)=0$. But also $y+b \in a_{s_2+1} \matn$.
\end{prose}

\item Recapitulating the choices so far,
\plist 
\item $\alpha(q)(s_1) > 0, \alpha (q) (s_2) < \frac{a_{s_2+1}}{a_{s_2}} -1$ (since $s_1 \in Q_1, s_2 \in Q_0 $).
\item $\alpha\xmax((\sum_{i\leq s_2} \alpha(q)(i) a_i) - a) >0$ (since $\alpha\xmax(a) < s_1\in Q_1$). \label{eq:positive}
\item $\alpha\xmax((\sum_{i\leq s_2} \alpha(q)(i) a_i) - a)$ is not maximal (since $\alpha(q)(s_2)$ not maximal and $s_2 > \alpha\xmax(a)$).
\item $\alpha\xmin(b+y) > s_2$ (since $b+y \in a_{s_2+1} \matn $).
\item $s_2 \notin \alpha\xsupp(a+y)$ (since $a + y  \in FS(\mathbf{x})  $).
\pliste 

\begin{prose}
The lurking contradiction lies in the fact that since $y$ translates such a small $a$ into $FS(\mathbf{x})$, it cannot simultaneously translate elements like $b$, i.e., elements that agree with $\alpha(q)$ up to $s_2$, to be divisible by $a_{s_2+1}$.

This is due to the (non-maximal) ``hole'' of both $(y+a)$ and $(b-a)$ at $s_2$ which simply does not allow for enough carrying over in the sum $(y+b)$ to get a multiple of $2^{s_2+1}$.
\end{prose}

\item First calculate
\begin{align*}
\sum_{i > s_2 } \alpha( b + y )(i)  a_i & = (a + y ) + (b - a)  \\
& = \sum_{i\in \omega } \alpha( a + y )(i) a_i + \sum_{i \in \omega } \alpha(b-a)(i) a_i , \\
\end{align*}
Recall that $b-a > 0$, so not all $\alpha(b-a)(i)$ are zero -- but $\alpha(b-a)(s_2)$ is not maximal (as noted before).

\item Rearranging this equation yields
\begin{align*}
& \sum_{i > s_2 } \alpha( b + y )(i)  a_i  - \sum_{i > s_2} \alpha( a + y )(i) a_i - \sum_{i > s_2} \alpha(b-a)(i) a_i \\
& = \sum_{i \leq s_2} \alpha( a + y )(i) a_i + \sum_{i \leq s_2} \alpha(b-a)(i) a_i  \\
& = \sum_{i < s_2} \alpha( a + y )(i) a_i +  \sum_{i \leq s_2} \alpha(b-a)(i) a_i ,
\end{align*}
since $s_2 \notin \alpha\xsupp(a+y)$.

\item Clearly, $a_{s_2+1}$ divides the first line, so the last line must add up to (a multiple of) $a_{s_2+1}$.

\begin{prose}
However, there is not enough carrying over.
\end{prose}
\item But
\[
0 < \sum_{i < s_2} \alpha( a + y )(i) a_i +  \sum_{i \leq s_2} \alpha(b-a)(i) a_i  < a_{s_2} + (\frac{a_{s_2+1}}{a_{s_2}}-1) a_{s_2} = a_{s_2+1}.
\]
\plist
\item Since $\alpha(b)$ agrees with $\alpha(q)$ up to $s_2$, step \ref{eq:positive} implies that both summands are positive.

\item Also since $\alpha(b)$ agrees with $\alpha(q)$ up to $s_2$, $\alpha(b)(s_2)$ is not maximal, i.e., less than $(\frac{a_{s_2+1}}{a_{s_2}}-1)$.

\item Finally, by choice of $s_1 > \alpha\xmax(a)$, also $\alpha(b-a)(s_2)$ is not maximal. 
\pliste
\item This contradiction completes the proof.
\pliste \end{proof} 

After this complicated proof, the main result follows almost immediately.

\begin{Theorem}[Trivial sums]\label{thm:trivial sums}
A strongly summable ultrafilter with disjoint $\alpha$-support can only be written as a sum trivially.
\end{Theorem}
\begin{proof}  \plist
\item Assume that $p$ is a strongly summable ultrafilter with disjoint $\alpha$-support and $q,r \in \beta \matn$ with
\[
q+r = p.
\]

\item The above Lemma \ref{lem:nearly trivial sums} implies $r \in \matz + \bigcap_{n\in \omega} \bbar{a_n\matn}$.

\item Therefore there exists $k\in \matz$ such that $-k+r \in \bigcap_{n\in \omega} \bbar{a_n\matn}$; in particular
\[
(k+q) + (-k +r ) = p.
\]

\item But now applying Lemma \ref{lem:trivial sums for H} with $k+q$ and $-k+r$ implies $k+q = -k +r =p$ -- as desired.
\pliste \end{proof} 

This result, however, leaves some obvious questions open.

\begin{Question}
\begin{itemize}
\item Does every strongly summable ultrafilter have the trivial sums property?
\item Does every strongly summable ultrafilter have disjoint $\alpha$-support for some $\mathbf{a}$?
\item Do other (idempotent) ultrafilters have the trivial sums property?
\end{itemize}
\end{Question}

A slight progress on the first two is the following proposition.

\begin{Proposition}\index{ultrafilter!{stable ordered union}}\index{ultrafilter!selective}
Let $p$ be a strongly summable ultrafilter additively isomorphic to a stable ordered union ultrafilter. Then $p$ has disjoint binary support (hence trivial sums).
\end{Proposition}
\begin{spoiler}
Ordered unions guarantee ordered $\mathbf{x}$-support for appropriate $\mathbf{x}$. Since $FS(\mathbf{x})$ always contains elements with ordered binary support, stability ``enforces'' this throughout a condensation.
\end{spoiler}

\begin{proof} \plist
\item Consider an additive isomorphism $\phii$ defined on a suitable $FS(\mathbf{x}) \in p$ such that $\phii(p)$ is stable ordered union.

\item Consider the following set
\[
\{ (v,w) \in \phii[FS(\mathbf{x})]_<^2 \ |\ \bmax(\phii^{-1}(v)) < \bmin(\phii^{-1}(w)) \}.
\]

\item Since $\phii(p)$ is a stable ordered union ultrafilter, there exists ordered $FU(\mathbf{s}) \in \phii(p)$ such that $FU(\mathbf{s})_<^2$ is included or disjoint from the above set.

\item But $FU(\mathbf{s})_<^2$ cannot be disjoint.
\plist
\item For any $FU(\mathbf{s}) \in \phii(p)$ there is some $\mathbf{y} \sqsubseteq \mathbf{x}$ with $\phii^{-1}[FU(\mathbf{s})] = FS(\mathbf{y})$.

\item But for any $z \in FS(\mathbf{y})$ we can pick $z' \in FS(\mathbf{y}) \cap 2^{\bmax(z)} \matn (\in p)$.

\item Then the pair $(\phii(z),\phii(z'))$ is included in the above set.
\pliste
\item The homogeneous $FU(\mathbf{s})\in \phii(p)$ yields some $\phii^{-1}[FU(\mathbf{s})] = FS(\mathbf{y}) \in p$.

\item Since $\mathbf{s}$ is ordered, $y$ must have ordered, hence disjoint binary support.
\pliste \end{proof} 

So, as usual, the strongest notion of strongly summable ultrafilter has the desired trivial sums property. A negative answer to the first question would probably require the identification of a new kind of union ultrafilter.

The most natural answer to the second question would be to prove that $\bsupp$ maps strongly summables to union ultrafilters -- after all, its inverse map maps union ultrafilters to strongly summable ultrafilters.

For the closing remark, recall the following two notions. An ultrafilter in $\mathbb{N}$ is a \emph{$P$-point} if whenever we pick countably many of its elements $(A_n)_{n \in \omega}$, it includes a pseudo-intersection $B$, i.e., $A_n \setminus B$ is finite for all $n$. An ultrafilter is \emph{rapid} if for every unbounded function $f:\mathbb{N} \rightarrow \mathbb{N}$ it contains an element $B$ such that $|f^{-1}(n) \cap B|\leq n$ for all $n$. Since union ultrafilters map to rapid $P$-points under $\max$, the following might suggest a positive answer.

\begin{Proposition}
Let $p$ be strongly summable. Then $\bmax(p)$ is a rapid $P$-point.
\end{Proposition}

\begin{spoiler}
The proof is a modification of the proof of \cite[Theorem 2]{BlassHindman87}.
\end{spoiler}

\begin{proof} \plist
\item Pick a sequence $\mathbf{x}$ with sufficient growth and $FS(\mathbf{x}) \in p$.

\item Given $f\in \omega^\omega$ consider the set
\[
A = \{ a \in FS(\mathbf{x})  \ | \ f(\bmax(a)) \leq \min( \mathbf{x} \xsupp (a) ) \}.
\]
Then either $A$ or its complement is in $p$.

\item If $A \in p$ then $f$ is bounded (and therefore constant) on a set in $\bmax(p)$.
\plist
\item Pick $a \in A$ and $FS(\mathbf{y}) \subseteq (A \cap FS_{>\mathbf{x}\xmax(a)} (\mathbf{x}) \cap 2^a\matn) $  in $p$.

\item Then for $b \in FS(\mathbf{y})$ calculate
\begin{align*}
f(\bmax(b)) = f(\bmax(a + b)) & \leq \min( \mathbf{x}\xsupp (a+b) )  \\
& = \min ( \mathbf{x}\xsupp (a)).
\end{align*}
\item In other words, $f$ is bounded on $\bmax[ FS ( \mathbf{y} ) ] \in \bmax(p)$.
\pliste

\item If $\matn \setminus A \in p$, then $f$ has $|f^{-1}(n) | \leq n$ on a set in $\bmax(p)$.
\plist
\item Pick $\mathbf{y} \sqsubseteq \mathbf{x}$ with $FS(\mathbf{y}) \in p$, disjoint from $A$.

\item Therefore, each $z \in FS(\mathbf{y})$ with $n = f(\bmax (z))$ must have $n > x \xmin(z)$.

\item Since $\mathbf{y}$ has sufficient growth, $\bmax[FS(\mathbf{y}) ] = \bmax[ \mathbf{y} ]$.

\item Due to the disjoint $\mathbf{x}$-support of the $y_i$, there are at most $n$ indices $i$ such that $n > x\xmin(y_i)$.
\pliste
\item This completes the proof.
\pliste \end{proof}  

Thanks to the above proposition we might favor that all strongly summable ultrafilters have disjoint binary support. However, an answer remains elusive. It seems, however, that further progress on writing strongly summable ultrafilters as sums might lead to a better understanding of the phenomena in $\beta \mathbb{N}$ in general, just as it did with strongly right maximality.

\section*{Acknowledgments}

This work evolved out of the author's Ph.D.~thesis \cite{Krautzberger09} written under the supervision of Sabine Koppelberg at the Freie Universität Berlin and supported by the NaF\"oG grant of the state of Berlin. The author is also very grateful for the support of Andreas Blass especially during a visit to the University of Michigan, Ann Arbor, in the winter 2007/2008 with the support of the DAAD.

\bibliographystyle{alpha}

\end{document}